\documentclass[review]{elsarticle}

\usepackage{hyperref}
\usepackage{lineno}
\usepackage{amsfonts}
\usepackage{amsmath}
\usepackage{amsxtra}
\usepackage{stmaryrd}
\usepackage{mathabx}
\usepackage{pdflscape}

    \newtheorem{theorem}{Theorem}[section]
    \newtheorem{lemma}[theorem]{Lemma}
    
    \newtheorem{corollary}[theorem]{Corollary}

    \newenvironment{proof}[1][Proof]{\begin{trivlist}
    \item[\hskip \labelsep {\bfseries #1}]}{\end{trivlist}}
    \newenvironment{definition}[1][Definition]{\begin{trivlist}
    \item[\hskip \labelsep {\bfseries #1}]}{\end{trivlist}}
    
    \newenvironment{notation}[1][Notation]{\begin{trivlist}
    \item[\hskip \labelsep {\bfseries #1}]}{\end{trivlist}}
    \newenvironment{example}[1][Example]{\begin{trivlist}
    \item[\hskip \labelsep {\bfseries #1}]}{\end{trivlist}}

    \newcommand{\A}{\mathcal{A}}
    \newcommand{\C}{\mathcal{C}}
    \newcommand{\D}{\mathcal{D}}
    \newcommand{\R}{\mathcal{R}}

\begin{document}

\begin{frontmatter}

\title{Solving systems of equations on antichains for the computation of the ninth Dedekind Number}

\author[a1]{Patrick De Causmaecker}
\address[a1]{KU Leuven, Department of Computer Science, KULAK, CODeS\\patrick.decausmaecker@kuleuven.be}
\author[a2]{Lennart Van Hirtum}
\address[a2]{Paderborn Center for Parallel Computing, Paderborn University\\lennart.vanhirtum@gmail.com}

\begin{abstract}
Recently the ninth Dedekind Number ($D(9)$) was computed (\cite{JAKEL2023100006,10.1145/3674147}). In fact, the result of two independent computations were published nearly at the same time, in one of them the authors of the present paper were involved. $D(n)$ counts the monotone Boolean functions or antichains on subsets of a set of $n$ elements. The number rises doubly exponentially in the number of elements $n$, and until now no algorithm of a lower combinatorial complexity is known to compute $D(n)$. In our computation, we use coefficients representing the number of solutions of a specific set of equations on antichains over a finite set. We refer to these coefficients as P-coefficients. These can be computed efficiently. In this paper, we generalise this coefficient and apply it to four different systems of equations. Finally we show how the coefficient was used in our computation of $D(9)$, and how its generalisations can be used to compute $D(n)$.
\end{abstract}

\begin{keyword}
Dedekind numbers, combinatorics, counting, monotone boolean functions, antichains, equations on antichains, algebra of distributive lattices, posets, intervals
\end{keyword}

\end{frontmatter}


\section{Introduction}
\label{sec:introduction}
In the last section of his seminal paper of 1897 \cite{Dedekind1897}, Richard Dedekind stated a problem of counting the number of elements in a particular partially ordered set. The sets were characterised by a number of generators, and he succeeded to count the number of elements for 3 generators (18) and for 4 generators (166). The problem is equivalent to counting the number of monotone boolean functions on subsets of finite sets of $n$ elements.  Later in the twentieth century, the numbers became known as the Dedekind numbers ($D(n)$). Counting the two trivial monotone boolean functions returning {\it false}, respectively {\it true}, on all subsets, which Dedekind initially did not, the two numbers from the 1897 paper are now known as the third (20) and the fourth (168) Dedekind number ($D(3)$ and $D(4)$). From then on, calculating Dedekind Numbers (OEIS series A000372 \cite{OEISDedekind}) has always followed progress in both the Mathematical domain, and the Computational Implementation Domain. In 1940 Randolph Church computed $D(5) = 7579$ (today 7581) \cite{CHURCH1940} using a new tabulation method. In 1946, Morgan Ward succeeded computing $D(6) = 7.828.352$ (today 7.828.354) \cite{WARD}. In 1965, Randolph Church used the CDC 1604 computer to arrive at $D(7) = 2.414.682.040.998$ \cite{CHURCH1965} after which it took 10 years until in 1976, Berman \& K\"ohler confirmed this using the IBM 370 \cite{BERMANKOHLER}. Since the year 1991 and up until March 2023, the highest known term in the sequence was $D(8) = 56.130.437.228.687.557.907.788$, computed by Doug Wiedemann \cite{wiedemannDedekind8} on a Cray 2 supercomputer. 

In April 2023 two independent computations of $D(9)$ were announced \cite{JAKEL2023100006, 10.1145/3674147}. The computations happened completely independently, were based on different methodologies and relied on very different hardware.  Christian J\"akel (\cite{JAKEL2023100006}) used a GPU infrastructure and built on a  matrix formulation, which allowed to use advanced techniques and highly advanced software for numerical analysis. In our computation (Van Hirtum et al. \cite{10.1145/3674147}), we, together with our co-authors, relied on the advanced, partially experimental, super computing infrastructure at the PC2 centre at Paderborn University, using cutting edge FPGA hardware and an ingenious implementation of the P-coefficient formula \cite{decausmaecker2014number}. At the end of the latter computation, some uncertainty due to the semi-experimental setting of the FPGA infrastructure was identified. Although we obtained the result on March 8, it was only when the result of J\"akel came out, confirming our result, that we could be certain. An extensive report on this computational achievement is in \cite{10.1145/3674147}, concentrating on the very important computational aspects and the innovative implementation which brought new insights for FPGA technology. The present paper elaborates on the mathematical aspects, with a discussion on the P-coefficient and possible extensions.

The structure of the paper is as follows. 
In Section \ref{sec:definitions}, we describe the formalism in which our theorems will be stated.
In Section \ref{sec:preliminaries}, we position our method in line with the approach used by Wiedemann in his computation of the eighth Dedekind number denoted as $D(8)$ \cite{wiedemannDedekind8}. 
In Section \ref{sec:related}, we describe some related work. 
In Section \ref{sec:twoequations}, we derive the P-coefficient result as we used it in our computation. In particular, we show how to count the number of solutions of a particular set of two equations with antichains as variables.
In Section \ref{sec:twoequationsapp}, we show how the result of Section \ref{sec:twoequations} can be used to compute $D(n+2)$ which led to our computation of $D(9)$.
In Section \ref{sec:general}, we generalise the result of Section \ref{sec:twoequations} on counting the number of solutions of a particular set of $k+1$ equations for any integer $k$.
In Sections \ref{sec:fourequations} and \ref{sec:fourequationsapp}, we apply the results of Section \ref{sec:general} for the computation of $D(n+3)$ using resuts on $D(n)$.
In Sections \ref{sec:sevenequations} and \ref{sec:sevenequationsapp}, we apply our main result to derive a formula for $D(n+4)$.
Finally, in Section \ref{sec:conclusions}, we briefly summarise the paper.

\section{Definitions and notations}
\label{sec:definitions}
\subsection{The Dedekind number}
For any finite positive number $n \ge 0$, the n-th Dedekind number counts the {\it number of monotone Boolean functions} (see paragraph \ref{def:monotonefunction}), or equivalently, the number of {\it antichains}, on the subsets of a finite set $\A = \{1,\dots,n\}$ (see paragraph \ref{def:antichain}).
To derive the properties in this paper, we will use the concept of an antichain over $\A$, i.e. sets of pairwise non-inclusive subsets of $\A$.
The set of all antichains over $\A$ will be denoted by $\D_n$.
The number of antichains over $\A$ is the $n^{th}$ Dedekind number and it will be denoted by $D(n)$.

The set of permutations of the elements of base set $\A$ generates an equivalence relation on $\D_n$.
The corresponding set of equivalence classes is denoted by $\R_n$ and the number of such equivalence classes is denoted by $R(n)$.

The known Dedekind numbers are shown in Table \ref{tab:dedekindNumbers}, including the result of our computation.
Table \ref{tab:equivalenceClassCounts} shows the known numbers $R(n)$ of equivalence classes of monotone Boolean functions under permutation of the elements of the base set.

\begin{table}[h]
    \centering
    \resizebox{\columnwidth}{!}{
    \begin{tabular}{c|l|c}
        D(0) & 2 & Dedekind (1897) \\
        D(1) & 3 & Dedekind (1897) \\
        D(2) & 6 & Dedekind (1897) \\
        D(3) & 20 & Dedekind (1897)\\
        D(4) & 168 & Dedekind (1897) \\
        D(5) & 7581 & Church (1940) \\
        D(6) & 7828354 & Ward (1946) \\
        D(7) & 2414682040998 & Church (1965) \\
        D(8) & 56130437228687557907788 & Wiedemann (1991) \\
        D(9) & 286386577668298411128469151667598498812366 & J\"akel, Van Hirtum et al. (2023)
    \end{tabular}
    }
    \caption{Known Dedekind Numbers, OEIS series A000372 \cite{OEISDedekind}, including our result.}
    \label{tab:dedekindNumbers}
\end{table}

\begin{table}[h]
    \centering
    \resizebox{\columnwidth}{!}{
    \begin{tabular}{c|l|c}
        R(0) & 2 & \\
        R(1) & 3 & \\
        R(2) & 5 & \\
        R(3) & 10 & \\
        R(4) & 30 & \\
        R(5) & 210 & \\
        R(6) & 16353 & \\
        R(7) & 490013148 & Tamon Stephen \& Timothy Yusun (2014) \cite{STEPHEN201415} \\
        R(8) & 1392195548889993358 & Bartelomiej Pawelski (2021) \cite{pawelski2021numberinequivalentmonotoneboolean}\\
        R(9) & 789204635842035040527740846300252680 & Bartelomiej Pawelski (2023) \cite{pawelski2023numberinequivalentmonotoneboolean}
    \end{tabular}
    }
    \caption{Known Equivalence Class Counts, OEIS series A003182 \cite{OEISDedekindNEQ}}
    \label{tab:equivalenceClassCounts}
\end{table}

\subsection{Notation}

Let $n \in \mathbb{N}, \A = \{1,\dots,n\}, 2^{\A} = \{X \subseteq \A\},\mathbb{B} = \{false,true\}$. A function $f:2^A \rightarrow \mathbb{B}$ is called a Boolean function on $2^A$.
\begin{definition}
\label{def:monotonefunction}
A Boolean function on $2^{\A}$ is monotone iff
$$\forall X,Y \in 2^{\A}, X \subset Y: f(Y) \Rightarrow f(X)$$
\end{definition}

\begin{definition}
\label{def:maximalset}
A set $X \in 2^{\A}$ is a maximal set of a monotone Boolean function $f$ iff
$$f(X) = true\ and\ \forall Y \in 2^{\A}: X \subsetneq Y \Rightarrow f(Y) = false$$
\end{definition}

A monotone Boolean function is completely defined by its set of maximal sets.

For any monotone Boolean function, no two of its maximal sets include each other.
A set of sets with this property is called an {\it antichain}:
\begin{definition}
\label{def:antichain}
A set of sets $\sigma \subseteq 2^{\A}$ is an antichain iff
$$\forall X,Y \in \sigma:X \not= Y \Rightarrow X \not\subseteq Y\ and\ Y \not\subseteq X$$
\end{definition}
In what follows, we will represent antichains by letters from the Greek alphabet $\alpha,\beta,\ldots$.
The elements of the set $\A$ will be denoted by Latin minuscules $a,b,c,\ldots$ or by digits $1 \ldots 9$. 
For subsets, unless there is a possibility for confusion, we will use the abbreviation $abc \equiv \{a,b,c\}, 123 \equiv \{1,2,3\}, 0 \equiv \emptyset$ or refer to them by Latin majuscules $X,Y,\ldots$. 

\begin{definition} An antichain $\alpha$ is said to dominate set $X \in 2^A$ if there is at least one set $Y \in \alpha$ such that $X \subseteq Y$.
\end{definition}
\begin{definition} An antichain $\alpha$ is less than or equal to an antichain $\beta$ iff each set in $\alpha$ is dominated by $\beta$:
$$\alpha \le \beta \Leftrightarrow \forall X \in \alpha: \exists Y \in \beta: X \subseteq Y$$
\end{definition}

\begin{notation}
We denote by $\bot$ and $\top$ the smallest, respectively the largest, element of $\D_n$: 
\begin{gather}
\bot \equiv \{\}, \top \equiv \{12\ldots n\}\\
\forall \alpha \in \D_n: \bot \le \alpha \le \top
\end{gather}
\end{notation}

\begin{definition}
For any two elements $\alpha,\beta \in \D_n$, the interval with bottom $\alpha$ and top $\beta$ is the set of antichains $\chi$ satisfying $\alpha \le \chi \le \beta$:
$$\forall \alpha, \beta \in \D_n: [\alpha,\beta] = \{\chi \in \D_n: \alpha \le \chi \le \beta\}$$
\end{definition}
\begin{definition}
For $\alpha,\beta \in \D_n$, the {\it join} $\alpha \vee \beta$ and the {\it meet} $\alpha \wedge \beta$ are the antichains defined by
\begin{align}
  \alpha \vee \beta &= max (\{X \in \alpha \cup \beta\})\\
  \alpha \wedge \beta &= max (\{X \cap Y | X \in \alpha, Y \in \beta\}) 
\end{align}
\end{definition}
Where $max$ denotes the maximum w.r.t. inclusion of a set of subsets of $\A$, i.e.
$$ S \subseteq 2^{\A} \Rightarrow max(S) = \{X \in S| \forall Y \in S:X \not\subsetneq Y\}$$

As an antichain is a set of sets, we use the operations intersection $(\cap)$, union $(\cup)$, subset $(\subset, \subseteq, \subsetneq)$ and difference $(-)$ to work on antichains.
These should not be confused with the join $(\vee)$, meet $(\wedge)$, and less or equal $(\le)$ operations in the partial order we have just introduced.
The span of an antichain  is the union of its elements:

\begin{definition}
\label{def:span}
$\forall \alpha \in \D_n:sp(\alpha) = \cup_{X \in \alpha} X$
\end{definition}

A last operator on two antichains is the direct product, denoted by '$\times$', and defined for any two antichains $\alpha$ and $\beta$ for which the span does not overlap:

\begin{definition}
\label{def:times}
$\forall \alpha,\beta \in \D_n, sp(\alpha) \cap sp(\beta) = \emptyset:\alpha \times \beta = \{X \cup Y | X \in \alpha, Y \in \beta\}$
\end{definition}

Finally, in the formulas below, a number defined for each pair $\alpha \le \beta \in \D_n$ plays an important role.
We refer to this number as the {\it connector number}  $C_{\alpha,\beta}$. It counts the number of connected components of a graph with  sets in $\beta - \alpha$ as vertices and where two sets are connected by an edge if the intersection is not dominated by $\alpha$ (exact definition follows).

\section{Preliminaries}
\label{sec:preliminaries}

In this section, we first describe the decomposition used by by Wiedemann \cite{wiedemannDedekind8} in his 1991 computation of $D(8)$.
Consequently we describe our own approach for the computation of $D(9)$ which was based on the same decomposition, except that we reformulated the summation formula and introduced a new counting factor. 

Wiedemann's decomposition was stated in the set of monotone boolean functions. In Wiedemann's notation, $Q(n)$ stands for the set of all subsets of $N = \{1,\ldots,n\}$ and a monotone boolean function is a subset $S$ of $Q(n)$ with the property $s \in S \Rightarrow \forall t \subset s: t \in S$. In this section, we will call such a set {\it monotonic}. He reduced the computation of $D(8)$ to a sum over sets in $D(6)$ by splitting a monotonic set $S$ in four parts as follows:
\begin{align*}
S_{00} &= \{s \in Q(6) : s \in S\}\\
S_{01} &= \{s \in Q(6) : (s \cup \{7\}) \in S\}\\
S_{10} &= \{s \in Q(6) : (s \cup \{6\}) \in S\}\\
S_{11} &= \{s \in Q(6) : (s \cup \{6,7\}) \in S\}
\end{align*}
Note that these four sets uniquely determine $S$, and are all monotonic if $S$ is monotonic.
Moreover 
\begin{gather}
\label{eq:twoinclusion}
S_{00} \subseteq S_{01},S_{00} \subseteq S_{10}, S_{01} \subseteq S_{11}, S_{10} \subseteq S_{11}
\end{gather}
The dual of $S$ can be defined as the complement of the set of complements of sets in S, or, again in Wiedemann's notation, $S^* = \{s^c | s \in S\}^c$. 
An integer valued function $\eta$ on monotonic sets is introduced, with $\eta(S)$ the number of monotonic subsets of a monotonic set $S$.
Wiedemann's basic formula now reads
\begin{gather}
\label{eq:wiedemann}
D(8) = \sum_{S_{01}\in \D_6}\sum_{S_{10} \in \D_6} \eta(S_{01} \cap S_{10})\eta(S^*_{01} \cap S^*_{10})
\end{gather}
where $\D_n$ stands for the set of all monotonic subsets over $Q(n)$.
The factor $\eta(S_{01} \cap S_{10})$ counts the number of sets $S_{00}$ allowed by $S_{01}$ and $S_{10}$ while the factor $\eta(S^*_{01} \cap S^*_{10})$ counts the number of allowed sets $S_{11}$.
Reducing one of the sums to a sum over  nonequivalent sets under permutation of the numbers $\{1,\ldots,6\}$, so over $R \in \R_6$, and introducing the number of equivalent sets in the corresponding equivalence class as $\gamma(R)$, the expression can be written as
\begin{gather*}
\label{eq:wiedemanneq}
D(8) = \sum_{R \in \R_6}\sum_{T \in \D_6} \gamma(R) \eta(R \cap T)\eta(R^* \cap T^*)
\end{gather*}
Berman and K\"ohler \cite{BERMANKOHLER} and Berman, Burger and K\"ohler \cite{BURGER}, develop a general binary-tree decomposition algorithm for counting the number of elements in free distributive lattices over partially ordered sets. They apply it to the case of a free distributive lattice on $n$ generators which is identical to the antichains over the set of subsets of a finite set of $n$ elements and hence to counting the Dedekind number. They report on the computation of $D(6)$ on an IBM 370, taking five minutes. To compute $D(7)$, they generate the elements of $D(5)$ using their decomposition algorithm and then apply exactly the same formula which was later used by Wiedemann \cite{wiedemannDedekind8} to compute $D(8)$ as explained above (\ref{eq:wiedemann}). They confirm the result obtained by Church in 1965 \cite{CHURCH1965} hereby closing the discussion raised by a conflicting result of a 1971 computation by Lunnon \cite{LUNNON}.

In our computation of $\D_9$ \cite{10.1145/3674147}, we reorder the terms in equation \ref{eq:wiedemann} and identify large numbers of identical terms. 
Counting these terms boils down to counting connected components of a specific set of graphs while the number of terms is reduced from $\approx 6\times 10^{13})$ to $\approx 2 \times 10^{12})$ in the case of $D(8)$ and from $\approx 2 \times 10^24$ to $\approx 5 \times 10^22$ in the case of $D(9)$. Moreover, our method allows permutation symmetry in the numbers $\{1,\ldots,9\}$ to be taken into account more effectively.

In what follows, we describe the methodology in detail.  
We start with the approach and formula's we used to compute $D(9)$. 
Consequently, we state and prove a generalisation identifying a number of possible approaches.
Finally, we apply the generalisation in specific cases, in particular to the formula's used by Christian J\"akel \cite{JAKEL2023100006}.

\section{Related Work}
\label{sec:related}
Very important to mention is the work of Christian J\"akel, who independently computed the 9th Dedekind Number at about the same time as us, but using different formulas and different hardware. His computation is based on the decomposition to be presented in Section \ref{sec:sevenequationsapp}. He did not use  P-coefficients, but his transformation of the formulas made it possible to use  matrix multiplication techniques which played an important role in this computation. This way he was able to compute D(9) in only 8000 GPU hours \cite{JAKEL2023100006}.

Bartłomiej Pawelski recently published some important work regarding a related sequence, namely the "Number of Inequivalent Monotone Boolean Functions". In 2021 he computed R(8) \cite{pawelski2021numberinequivalentmonotoneboolean}, quickly followed in 2023 by R(9), which was dependent on our and J\"akel's D(9) finding. \cite{pawelski2023numberinequivalentmonotoneboolean}. Afterwards, he has worked on further related sequences, such as the number of self-dual Monotone Boolean Functions \cite{pawelski2023countingselfdualmonotoneboolean}. 
Another result of Bartolomiej Pawelski is on properties of the Dedekind numbers, see \cite{pawelski2023divisibility}.

The work in the present article builds and extends upon three papers by Patrick De Causmaecker et al. \cite{decausmaecker2014number,decausmaecker2011partitioning,decausmaecker2016intervals}.

\section{A system of two equations}
\label{sec:twoequations}
 \label{sec:meet2join}
Let $\alpha$ and $\beta$ be two antichains with $\alpha \le \beta$.
Consider the set of equations, hereafter refered to as System II.
\begin{align}
 \chi \wedge \upsilon &= \alpha \label{meet2constraint}\\
 \chi \vee \upsilon &= \beta \label{join2constraint} 
\end{align}
It is clear that the system has no solution if $\alpha \not\le \beta$.
The following Lemma states that in any solution $\chi,\upsilon$, only the sets in the righthand sides of System II can occur.
This lemma significantly simplifies the construction as well as the counting of solutions of System  II .
\begin{lemma}
\label{lem:onlysetsofbeta}
In any solution of System II, $\chi = \alpha \vee \chi_1, \upsilon = \alpha \vee \upsilon_1$ with $\chi_1 \subseteq \beta-\alpha, \upsilon_1 \subseteq \beta-\alpha$.
\begin{proof}
Let $\chi_1 = \chi - \alpha$. Suppose $X \in \chi_1, X \not\in \beta$. Due to equation \ref{join2constraint}, there must be a set $X' \in \upsilon$ such that $X \subseteq X'$ and $X' \in \beta$. Due to the equation \ref{meet2constraint}, this implies $X \in \alpha$, a contradiction. Similar arguments work for $\upsilon_1$. \qed
\end{proof}
\end{lemma}
The following lemma states that any set satisfying the condition of Lemma \ref{lem:onlysetsofbeta} occurs in exactly one of the two antichain variables $\chi$ or $\upsilon$.
\begin{lemma}
\label{oneOfTwo}
For any set $A \in \beta - \alpha$, we have that $A \in \chi - \upsilon$ or $A \in \upsilon - \chi$.
\begin{proof}
From equation \ref{join2constraint}, it follows that $A$ must be in at least one of $\chi$ or $\upsilon$.
From the same equation it follows that no set dominating $A$ can be in any of the two.
If $A$ were in both $\chi$ and $\upsilon$, $A \in \chi \wedge \upsilon$ and equation \ref{meet2constraint} would be violated.
\qed
\end{proof}
\begin{lemma}
\label{connect2}
Let $A,B \in \beta - \alpha$, such that $\{A \cap B\} \not\le \alpha$. 
We have
\begin{align*}
A \in \chi &\Leftrightarrow B \in \chi\\ 
A \in \upsilon &\Leftrightarrow B \in \upsilon
\end{align*}
\begin{proof}
Due to Lemma \ref{oneOfTwo}, $A$ and $B$ are each in exactly one of $\chi$ and $\upsilon$. If they were each in a different one, we would have $ \{A \cap B\} \le \chi \wedge \upsilon = \alpha$ due to equation \ref{meet2constraint}, a contradiction.
\qed 
\end{proof}
\end{lemma}
\end{lemma}
\begin{definition}
Given two antichains $\alpha \le \beta$, two sets $A,B \in \beta - \alpha$ are said to be {\it directly connected} with respect to $\alpha$ iff $\{A \cap B\} \not\le \alpha$. 
\end{definition}
\begin{definition}
Two sets $A,B \in \beta - \alpha$ are said to be connected w.r.t. $\alpha$ iff $A$ and $B$ are directly connected w.r.t. $\alpha$ or there exist $Z_1,Z_2,\ldots,Z_x \in \beta$, $x \ge 1$, such that $A$ is directly connected with $Z_1$, $B$ is directly connected with $Z_x$ and each $Z_i$ is directly connected to $Z_{i+1}$ for $i \in \{1,\ldots,x-1\}$, all w.r.t. $\alpha$.
\end{definition}
\begin{definition}
\label{connectiongraph}
For two antichains $\alpha \le \beta$, the connection graph $C_{\alpha}(\beta)$ is the graph with the sets of $\beta-\alpha$ as vertices and edges $$E = \{(A,B) | A,B \in \beta-\alpha, A,B \text{ are connected w.r.t. }\alpha\}$$
\end{definition}
\begin{theorem}
\label{pcoeffTheorem}
Let $\alpha \le \beta$ be two antichains, $(\chi_0,\upsilon_0)$ a partition of $\beta-\alpha$ such that any two sets $A,B \in \beta-\alpha$ that are connected w.r.t. $\alpha$ are in the same component of the partition. With the partition $(\chi_0,\upsilon_0)$ corresponds exactly one solution of System II and any solution of System II corresponds with exactly one such partition.
\begin{proof}
Due to Lemma \ref{lem:onlysetsofbeta}, any solution is determined by two subsets of $\beta - \alpha$.
Due to lemma \ref{connect2}, any two directly connected sets $A,B \in \beta$ must come in the same antichain $\chi$ or $\upsilon$. By transitive closure, this must be true for any two connected sets $A,B \in \beta$ and hence for all sets in any connected component of $C_{\alpha}(\beta)$. Given a partition $(\chi_0,\upsilon_0)$  of $\beta - \alpha$ satisfying the criteria of the theorem, we arrive at a solution $\chi = \chi_0 \vee \alpha, \upsilon = \upsilon_0 \vee \alpha$. The reverse follows by identifying $\chi_0 = \chi - \alpha, \upsilon_0 = \upsilon - \alpha$. \qed
\end{proof}
\end{theorem}
\begin{corollary}
The number of solutions of System II is given by $$P(\alpha,\beta) = 2^{C_{\alpha,\beta}}$$ where $C_{\alpha,\beta}$ is the number of connected components of the connection graph $C_{\alpha}(\beta)$.
\begin{proof}
The number of partitions in two sets of the connected components of $C_{\alpha}(\beta)$ is given by $2^{C_{\alpha,\beta}}$. As was proven in Theorem \ref{pcoeffTheorem}, this is the number of solutions of System II.
\end{proof}
\end{corollary}
\begin{notation}
The power of $2$, $2^{C_{\alpha,\beta}}$, is referred to as the {\it P-coefficient} and is denoted by $P(\alpha,\beta)$.
\end{notation}
\begin{corollary}
The number of solutions of System II can be computed in polynomial time in the number of sets in $\beta - \alpha$.
\begin{proof}
By Theorem \ref{pcoeffTheorem}, the problem is reduced to counting the connected components in an undirected graph with sets in $\beta - \alpha$ as vertices. This can be done in $O(V+E)$ time where $V$ is the number of vertices and $E$ is the number of edges. \qed
\end{proof}
\end{corollary}

\section{Application of System II }
\label{sec:twoequationsapp}
System II in Section \ref{sec:meet2join} can be used to compute the $(n+2)^{nd}$ Dedekind number $D(n+2)$ based on intervals in the space of antichains over $n$ elements.
\begin{theorem}
\label{the:nplus2}
Let $n \ge 0$ be a natural number, $N=\{1,..,N\}$, $\A_n$ the lattice of antichains over $N$. We have
$$ D(n+2) = \sum_{\alpha \le \beta \in \D_n}
2^{C_{\alpha,\beta}}  |[\bot,\alpha]||[\beta,\top]|$$
\begin{proof}
We use the decomposition for $\eta \in \D_{n+2}$ in terms of elements $\eta_X$ of $\D_n$:
\begin{align*}
 \eta  = & \eta_0 \times \{12\}\\
 & \vee  \eta_1 \times \{1\} \vee \eta_2 \times \{2\}\\
 & \vee  \eta_{12} \times \{0\}
 \end{align*}
where sets such as $\{1,2\}$ are denoted by $12$ and $\emptyset$ by $0$.
The decomposition is unique iff $A \subseteq B \Rightarrow \eta_A \le \eta_B$.
We identify $\eta_1,\eta_2$ with the variables $\chi,\upsilon$ in System II.
Any solution of System II for specific values of $\alpha,\beta$ will now correspond
to exactly one antichain $\eta$ for which $\eta_0 \in [\bot,\alpha]$, $\eta_{12} \in [\beta,\top]$. \qed
\end{proof}
\end{theorem}

\begin{corollary}
\label{cor:nplus2}
$D(n+2)$ can be computed as a sum of $D(n+1)$ terms, each of which is the product of the sizes of two intervals in $\D_n$.
\begin{proof} The condition $\alpha \le \beta$ for $\alpha,\beta \in \D_n$ allows exactly $D(n+1)$ terms.\qed
\end{proof}
\end{corollary}

This reduces the number of terms used in the computation of $D(8)$ by Wiedemann \cite{wiedemannDedekind8}, because without the condition $\alpha \le \beta$, the number of terms is $D(n)^2 $ ($D(6)^2 \approx 6\times 10^{13}$ while  $D(7) \approx 2 \times 10^{12}$ (Table \ref{tab:dedekindNumbers})).
Theorem \ref{the:nplus2} and its corollary \ref{cor:nplus2} were the basis for our computation of $D(9)$\cite{10.1145/3674147}. In this case we have a reduction in the number of terms from $D(7)^2 \approx 4 \times 10^{24}$ to $D(8) \approx 6 \times 10^{22}$.
Of course in all computations, symmetry is essential. The sum in Theorem \ref{the:nplus2} is essentially over intervals in $D(n)$ which on its own already allows to reduce the number of terms further to $\approx 10^{18}$ ($R(8)$ in table \ref{tab:equivalenceClassCounts}). More details can be found in Christian J\"akel's and our detailed reports on the computation of the $D(9)$ \cite{JAKEL2023100006, 10.1145/3674147}.

\section{General case}
\label{sec:general}
\label{sec:meetrjoin}
Let $\alpha,\{\beta_{ij}|1 \le i \le r,1 \le j \le r,  i \not= j\}\}$ be a set of antichains for some integer $r > 1$.
We will assume $\beta_{ij} = \beta_{ji}$ and will not distinguish between them.
We consider the following set of $1 + \frac{r(r+1)}{2}$ equations in $r$ variables $\chi_k, 1 \le k \le r$, hereafter referred to as System G :
\begin{align}
\bigwedge_{i=1}^{r} \chi_i & = \alpha           \label{meetrconstraint}\\
\chi_k \vee \chi_l & = \beta_{kl}\ \ \ \ \  (1 \le k \le r,1 \le l \le r, k \not= l) \label{joinrconstraint}
\end{align}

\begin{notation}
We will use $\beta = \cup_{k \not= l} \beta_{kl}$, $\gamma = \beta - \alpha$.
\end{notation}

Note that for System G to have a solution, no set in $\gamma$ can be dominated by $\alpha$.
Given equation \ref{joinrconstraint}, for any solution of System  G, for $k \not= l \in \{1,\ldots,r\}$, each set in $\beta_{kl}$ must be in $\chi_k$ or $\chi_l$. 

The following Lemma proves the reverse is also true: the only sets occurring in a solution of System  G are sets in $\beta$.
\begin{lemma}
\label{lem:onelybetasets}
In any solution of System  G,  $\chi_k \subseteq \beta$ for $k \in \{1,\ldots,r\}$.
\begin{proof}
Let $X \in \chi_k, X \not\in \gamma$. 
According to the equation \ref{joinrconstraint} $\chi_k \vee \chi_i = \beta_{ki}$ for any $i \not= k$ and there must be a set $Y_i \in \beta_{ki}$ s.t. $X \subsetneq Y_i$.
To satisfy the equation, we must have $Y_i \in \chi_i$.
According to equation \ref{meetrconstraint}, we find $\alpha = \chi_k \wedge (\wedge_{i \not= k} \chi_i) \ge \{X\}$ and hence $X \in \alpha$. \qed
\end{proof}
\end{lemma}

\begin{definition}
For any antichain $\delta$, a set $X$ is said to be dominated by $\delta$ if $\{X\} \le \delta$.
\end{definition}
\begin{definition}
A set $X$ is said to be dominating in a set of sets $\rho$ if $X \in \rho$ and there is no set $Y \in \rho$ s.t. $X \subsetneq Y$.
\end{definition}


\begin{lemma}
\label{lem:noteverywhere}
Let a set $X$ be dominating in $\gamma$ with at least one $\beta_{kl}$ such that $X \not\in \beta_{kl}$. If the System  G has at least one solution , there is a set $S_X \subseteq \{1,\ldots,r\}, S_X \not= \emptyset$ such that 
\begin{align*}
s \in S_X & \Leftrightarrow \forall i \in \{1,\ldots,r\} - \{s\} : X \in \beta_{si}
\end{align*}
and, for any solution $\{\chi_i| 1 \le i \le r\}$
\begin{align*}
s \in S_X & \Leftrightarrow \forall i \in \{1,\ldots,r\} - \{s\} : X \in \chi_s
\end{align*}
\begin{proof}
$X \not\in \beta_{kl} \Rightarrow X \not\in \chi_k\text{ and }X \not\in \chi_l$.
$X$ is in at least one of $\beta_{uv}$ and hence is in $\chi_u$ or $\chi_v$. 
If $X \in \chi_u$ then $u \not= k,u \not= l$ and $X \in \beta_{ui}$ for any $i \in \{1,\ldots,r\}$.
In particular $X \in \beta_{uk}$ and $X \in \beta_{ul}$ implying that $X$ must be in $\chi_u$ in any solution.
The same reasoning follows if $X \in \chi_v$ and the lemma follows. \qed
\end{proof}
\end{lemma}

\begin{lemma}
\label{lem:everywhere}
For any set $X$ dominating in $\gamma$, if $X \in \beta_{ij}$ for all $1 \le i < j \le r$, then in any solution of System G, there is exactly one $\chi_k$, $1 \le k \le r$, such that $X \not\in \chi_k$.
\begin{proof}
Suppose $X \not\in \chi_k, X \not\in \chi_l, k \not= l$. Then $X \not\in \beta_{kl}$, in contradiction with the condition in the lemma.
Suppose $X \in \chi_k$ for all $1 \le k \le r$. Then equation \ref{meetrconstraint} is violated. \qed
\end{proof}
\end{lemma}

\begin{definition}
\label{def:directlyconnected}
Two sets $X,Y \in \gamma$ are said to be directly connected w.r.t $\alpha$ if $\{X \cap Y\} \not\le \alpha$
\end{definition}

\begin{definition}
\label{def:connected}
Two sets $X,Y \in \gamma$ are said to be connected w.r.t $\alpha$ if they are directly connected w.r.t $\alpha$ or if there is a sequence $Z_1,Z_2,\ldots,Z_k$ of sets in $\gamma$
such that $X = Z_1, Y = Z_k$ and $Z_i,Z_{i+1}$ are directly connected for $1 \le i < k$.
\end{definition}

\begin{definition}
\label{def:graphofsets}
We denote by $g_{\alpha,\gamma}$ the graph with sets in $\gamma$ as vertices where there is an edge between two sets $X,Y$ iff $X$ and $Y$ are connected with respect to $\alpha$. 
\end{definition}

\begin{definition}
\label{def:connectedcomponents}
We denote by $C_{\alpha,\gamma}$ set of connected components of $g_{\alpha,\gamma}$
\end{definition}

The connected components of $g_{\alpha,\gamma}$ can be found in time of $O(V+E)$ where $V$ is the number of vertices in $g_{\alpha,\gamma}$, i.e. the number of sets in $\gamma - \alpha$ and $E$ the number of edges (\cite{Tarjan1972}).

\begin{definition}
\label{def:graphweight}
For a connected component $c \in C_{\alpha,\gamma}$, using Lemma \ref{lem:noteverywhere}, its weight $w(c)$ is defined as follows:
\begin{itemize}
\item[-] In case $c$ does not contain a set occurring in all $\beta_{ij}$: $$\text{if }\cup_{X \in c} S_X = \{1,\ldots,r\}\text{ then }w(c) = 0\text{ else }w(c) = 1$$
\item[-]In case $c$ contains at least one set dominated by all $\beta_{ij}$: $$w(c) = r - |\cup_{X \in c} S_X|$$
\end{itemize}
\end{definition}
 
We now prove the following
\begin{theorem}
\label{the:general}
The number of solutions of System G is
\begin{gather}
P_r (\alpha,\beta_{ij}\text{ for }1 \le i < j \le r) = \Pi_{c \in C_{\alpha,\gamma}} w(c)
\end{gather}
\begin{proof}
Note that, in any solution of System G for each connected component $c \in \C_{\alpha,\gamma}$, there must be at least one $\chi_i$
not containing any of the sets in $c$ because if this were not the case, equation \ref{meetrconstraint} would be violated.
We now reduce the system set by set as follows
\begin{enumerate}
\item Copy $\gamma$ to $\gamma'$, all $\beta_{ij}$ to $\beta'_{ij}$. Let all $\chi_i = \{\}$ be empty antichains.
\item \label{enu:starting} Choose a dominating set $X$ in $\gamma'$ which does not occur in all $\beta'_{ij}$. 
         If no such set exists goto step \ref{enu:solved}
\item Apply Lemma \ref{lem:noteverywhere} to assign the set to the right $\chi_i$'s.
\item Remove $X$ from $\gamma'$ and from each $\beta'_{ij}$.
\item For any strict subset $Y \subsetneq X$ present in $\gamma'$, replace each $\beta'_{ij}$ by $\beta'_{ij} \vee \{Y\}$.
\item Go to step \ref{enu:starting}
\item \label{enu:solved} Compute the weight of each $c \in C_{\alpha,\gamma}$ using the constructed $\chi_i$ as partial solutions.
\end{enumerate}
The formula in the theorem now counts the number of solutions. \qed
\end{proof}
\end{theorem}

\section{A system of four equations}
\label{sec:fourequations}
Let $\alpha \le \beta_{23},\beta_{13},\beta_{12}$ be four antichains with $\alpha \le \beta_i$ for $i \in {1,2,3}$.
Consider the set of equations, herafter referred to as System IV:
\begin{eqnarray}
 \chi_1 \wedge \chi_2 \wedge \chi_3 = \alpha \label{meet4constraint}\\
 \chi_1 \vee \chi_2 = \beta_{12} \label{join4constraint1} \\
 \chi_1 \vee \chi_3 = \beta_{13} \label{join4constraint2} \\
 \chi_2 \vee \chi_3 = \beta_{23} \label{join4constraint3} 
\end{eqnarray}

This is equivalent to the general System G from Section \ref{sec:meetrjoin} for $r = 3$.
The number of solutions is given by
$$P_3(\alpha,\beta_{12},\beta_{13},\beta_{23}) = \prod_{\substack{g\text{ connected component}\\ \text{ of }C_{\alpha}(\gamma)}}\ w(g)$$

\begin{corollary}
The number of solutions of System IV can be computed in polynomial time in the number of sets in $\beta - \alpha$.
\end{corollary}

\section{Application of System IV}
\label{sec:fourequationsapp}
System IV from Section \ref{sec:fourequations} can be used to compute the $(n+3)^{th}$ Dedekind number $D(n+3)$ based on intervals in the space of antichains over $n$ elements.
\begin{theorem}
Let $n \ge 0$ be a natural number, $\A=\{1,..,n\}$, $\D_n$ the lattice of antichains over $\A$. We have
$$ D(n+3) = \sum_{\substack{\alpha \in \D_n\\
\beta_{12},\beta_{13},\beta_{23} \in [\alpha,\top]\\
\gamma \in [\beta_{12} \vee \beta_{13} \vee \beta_{23},\top]}}
P_3(\alpha,\beta_{23},\beta_{13},\beta_{12})  |[\bot,\alpha]||[\beta_{12},\gamma]||\beta_{13},\gamma]||[\beta_{23},\gamma]$$
\begin{proof}
We use the decomposition for $\eta \in \D_{n+3}$ in terms of elements of $\D_n$:
\begin{align*}
 \eta  = & \eta_0 \times \{123\}\\
 & \vee  \eta_1 \times \{23\} \vee \eta_2 \times \{13\} \vee \eta_3 \times \{12\}\\
 & \vee  \eta_{12} \times \{3\} \vee \eta_{13} \times \{2\} \vee \eta_{23} \vee \{1\} \\
 & \vee  \eta_{123} \times \{0\}
 \end{align*}
were sets such as $\{1,3\}$ are denoted by $13$ and $\emptyset$ by $0$ as usual.
The decomposition is unique iff $A \subseteq B \Rightarrow \eta_A \le \eta_B$.
We identify $\eta_1,\eta_2,\eta_3$ with the variables $\chi_1,\chi_2,\chi_3$ in System IV.
Any solution of System IV for specific values of $\alpha,\beta_{12},\beta_{13},\beta_{23}$ will now correspond
to antichains $\eta$ for which $\eta_0 \le \alpha$, $\eta_{123} \ge \gamma$ for any $\gamma \ge \beta_{23} \vee \beta_{13} \vee \beta_{12}$ and $\eta_{12} \in [\beta_{12},\gamma],\eta_{13} \in [\beta_{13},\gamma],\eta_{23} \in [\beta_{23},\gamma]$. \qed
\end{proof}
\end{theorem}

\begin{example}
Consider only the antichains $\bot$ and $\{0\}$.
Table \ref{table:simplestP3} shows the resulting inequivalent instances of System IV, the number of equivalents under permutation and their number of solutions in the column $P_3$. The result is $D(3)$.
\end{example}
\begin{example}
Consider only the antichains in $[\bot,\{1\}]$.
Table \ref{table:obsP3} shows the resulting inequivalent instances of System IV, the number of equivalents under permutation and their number of solutions in the column $P_3$. The result is $D(4)$.
\end{example}

\begin{table}
   \center
   \resizebox{\columnwidth}{!}{
   \begin{tabular}{|c|c|c|c|c|c|c|c|}
   \hline
   $\alpha$ & $\beta_{12}$ & $\beta_{13}$ & $\beta_{23}$ & \#Eq & $P_3$ & 
   $\sum_{\gamma \in [\beta_{12} \vee \beta_{13} \vee \beta_{23},\{0\}]} |[\bot,\alpha]||[\beta_{12},\gamma]||[\beta_{13},\gamma]||[\beta_{23},\gamma]|$ & Tot \\
   \hline
   $\bot$         & $\bot$            & $\bot$            & $\bot$         & 1   & 1 & 1 + 8 & 9\\
   \hline
   $\bot$         & $\{0\}$          & $\{0\}$          & $\bot$         & 3   & 1 & 2 & 6\\
   \hline
   $\bot$         & $\{0\}$          & $\{0\}$          & $\{0\}$       & 1   & 3 & 1 & 3\\
   \hline
   $\{0\}$       & $\{0\}$          & $\{0\}$          & $\{0\}$       & 1  &  1 & 2 & 2\\
   \hline
   \multicolumn{7}{|c|}{$\sum $} & 20 \\
   \hline 
\end{tabular}
}
\caption{Simplest example of System \ref{sec:fourequations}, the result is D(3)}
\label{table:simplestP3}
\end{table}

\begin{table}
   \center
   \resizebox{\columnwidth}{!}{
   \begin{tabular}{|c|c|c|c|c|c|c|c|}
   \hline
   $\alpha$ & $\beta_{12}$ & $\beta_{13}$ & $\beta_{23}$ & \#Eq & $P_3$ & 
   $|[\bot,\alpha]|\times\sum_{\gamma \in [\beta_{12} \vee \beta_{13} \vee \beta_{23},\{1\}]} |[\beta_{12},\gamma]||[\beta_{13},\gamma]||[\beta_{23},\gamma]|$ & Tot\\
   \hline
   $\bot$         & $\bot$            & $\bot$            & $\bot$         & 1   & 1 & 1 + 8 + 27 = 36 & 36\\
   \hline
   $\bot$         & $\{0\}$          & $\{0\}$          & $\bot$         & 3   & 1 & 2 + 12 = 14 & 42\\
   \hline
   $\bot$         & $\{0\}$          & $\{0\}$          & $\{0\}$       & 1   & 3 & 1 + 8 = 9 & 27\\
   \hline
   $\bot$         & $\{1\}$          & $\{1\}$          & $\bot$         & 3  &  1 & 3 & 9\\
   \hline
   $\bot$         & $\{1\}$          & $\{1\}$          & $\{0\}$         & 3  &  2 & 2 & 12\\
   \hline
   $\bot$         & $\{1\}$          & $\{1\}$          & $\{1\}$         & 1  &  3 & 1 & 3\\
   \hline
   $\{0\}$         & $\{0\}$          & $\{0\}$          & $\{0\}$       & 1   & 1 & 2 $\times$ (1 + 8) = 18 & 18\\
   \hline
   $\{0\}$         & $\{1\}$          & $\{1\}$          & $\{0\}$         & 3  &  1 & 2$\times$2 = 4 & 12\\
   \hline
   $\{0\}$         & $\{1\}$          & $\{1\}$          & $\{1\}$         & 1  &  3 & 2$\times$1 = 2 & 6\\
   \hline
   $\{1\}$         & $\{1\}$          & $\{1\}$          & $\{1\}$         & 1  &  1 & 3$\times$1 = 3 & 3\\
   \hline
   \multicolumn{7}{|c|}{$\sum $} & 168 \\
   \hline 
\end{tabular}
}
\caption{One but simplest example of System \ref{sec:fourequations}, the result is D(4)}
\label{table:obsP3}
\end{table}

\section{A system of seven equations}
\label{sec:sevenequations}
Let $\alpha$ and $\beta_{ij}$ with $\alpha \le \beta_{ij}$ for $\{i,j\} \in \{1,2,3,4\}$ be seven antichains.
Note that we will use $ij$ to denote the combination of $i$ and $j$ without any order. We have e.g. $\beta_{ij} \equiv \beta_{ji}$ and $\beta_{1j} \equiv \beta_{j1}$.
Consider the set of seven equations in four variables $\chi_{1},\chi_{2},\chi_{3},\chi_{4}$, hereafter referred to as System VII:
\begin{eqnarray}
 \chi_{1} \wedge \chi_{2} \wedge \chi_{3} \wedge \chi_{4} = \alpha \label{meet7constraint1}\\
 \chi_{1} \vee \chi_{2} = \beta_{12} \label{join7constraint2}\\
 \chi_{1} \vee \chi_{3} = \beta_{13} \label{join7constraint3}\\
 \chi_{1} \vee \chi_{4} = \beta_{14} \label{join7constraint4}\\
 \chi_{2} \vee \chi_{3} = \beta_{23} \label{join7constraint5}\\
 \chi_{2} \vee \chi_{4} = \beta_{24} \label{join7constraint6}\\
 \chi_{3} \vee \chi_{4} = \beta_{34} \label{join7constraint7}
\end{eqnarray}

We will use the notation $\beta = \beta_{12} \cup  \beta_{13} \cup  \beta_{14} \cup  \beta_{23} \cup  \beta_{24} \cup \beta_{34}$ and $\gamma = \beta - \alpha$. System VII is a case of System G with $r = 4$. The number of solutions for System  \ref{sec:sevenequations} is given by
$$P_4(\alpha,\beta_{12},\beta_{13},\beta_{14},\beta_{23},\beta_{24},\beta_{34}) = \Pi_{\kappa \in C_{\alpha,\gamma}} w(\kappa)$$
\section{Application of System VII}
\label{sec:sevenequationsapp}
System VII from Section \ref{sec:sevenequations} can be used to compute the $(n+4)^{th}$ Dedekind number $D(n+4)$ based on intervals in the space of antichains over $n$ elements.
\begin{theorem}
\label{the:twotimesseven}
Let $n \ge 0$ be a natural number, $\A=\{1,..,n\}$, $\D_n$ the lattice of antichains over $\A$. We have
\begin{align*}
 D(n+4) =  & \sum_{\substack{
\alpha \le \beta_{12} \wedge \beta_{13} \wedge \beta_{14} \wedge \beta_{23} \wedge \beta_{24} \wedge \beta_{34}\\
\epsilon \ge \delta_{12} \vee \delta_{13} \vee \delta_{14} \vee \delta_{23} \vee \delta_{24} \vee \delta_{34}\\
\beta_{ij} \le \delta_{ij}
}} \nonumber \\
& P_4(\alpha,\beta_{12},\beta_{13},\beta_{14},\beta_{23},\beta_{24},\beta_{34}) \nonumber\\
\times & P_4(\widetilde{\epsilon},\widetilde{\delta_{12}},\widetilde{\delta_{13}},\widetilde{\delta_{14}},
\widetilde{\delta_{23}},\widetilde{\delta_{24}},\widetilde{\delta_{34}})\nonumber\\
\times & |[\bot,\alpha_]| \nonumber \\
\times&  |[\beta_{12},\delta_{12}]| \times |[\beta_{13},\delta_{13}]| \times |[\beta_{14},\delta_{14}]| \times |[\beta_{23},\delta_{23}]| \times |[\beta_{24},\delta_{24}]| \times |[\beta_{34},\delta_{34}]| \nonumber \\
\times & |[\epsilon,\top]|
\end{align*}
\begin{proof}
We use the decomposition for $\eta \in \D_{n+4}$ in terms of elements of $\D_n$:
\begin{align*}
 \eta  = & \eta_0 \times \{1234\}\\
 & \vee  \eta_1 \times \{234\} \vee \eta_2 \times \{134\} \vee \eta_3 \times \{124\} \vee \eta_4 \times {123}\\
& \vee  \eta_{12} \times \{34\} \vee \eta_{13} \times \{24\} \vee \eta_{14} \times \{23\} 
\vee \eta_{23} \times \{14\}  \vee \eta_{24} \times \{13\} \vee \eta_{34} \times \{12\}\\
 & \vee  \eta_{123} \times \{4\} \vee \eta_{124} \times \{3\} \vee \eta_{134} \times \{2\} \vee \eta_{234} \times \{1\} \\
 & \vee  \eta_{1234} \times \{0\}
 \end{align*}
The decomposition is unique iff $A \subseteq B \Rightarrow \eta_A \le \eta_B$.
In each term of the sum in the righthand side in the theorem, two instances of System VII are used.
In the one counted by 
$P_4(\alpha,\beta_{12},\beta_{13},\beta_{14},\beta_{23},\beta_{24},\beta_{34})$,
the variables $\chi_{i}$ are identified with $\eta_{1},\eta_{2},\eta_{3}, \eta_{4}$.
In the one counted by 
$P_4(\widetilde{\epsilon},\widetilde{\delta_{12}},\widetilde{\delta_{13}},\widetilde{\delta_{14}},\widetilde{\delta_{23}},\widetilde{\delta_{24}},\widetilde{\delta_{34}})$,
the variables $\chi_{i}$ are identified with the duals of $\eta_{123},\eta_{124},\eta_{134}, \eta_{234}$.
Any combined solution of the two realisations of System VII for specific values of $\alpha,\epsilon, \beta_{ij}$ and $\delta_{ij}$ can now be seen to correspond to antichains $\eta$ for which $\eta_0 \le \alpha$, $\eta_{1234} \ge \epsilon$, $\beta_{ij} \le \eta_{ij} \le \delta_{ij}$. 
The product of the two $P_4$ coefficients counts the number of combined solutions for $\eta_{i}$ and $\eta_{ijk}$.\qed
\end{proof}
\end{theorem}

Christian J\"akel (\cite{JAKEL2023100006}) in his computation used a technique based on the decomposition used in the proof of Theorem \ref{the:twotimesseven}. He  implemented it on his GPU infrastructure through matrix multiplication.
He did not use the P-coefficient formula proven here. 
The result in Theorem \ref{the:twotimesseven} requires to sum over fourteen variables which are constrained by a strong partial order.
This reduces the number of terms and allows the sum to be factorised.

\begin{example}
We apply Theorem \ref{the:twotimesseven} to its simplest case of computing $D(4)$. 
Table \ref{table:twosevensbeta} shows the non-equivalent values for $\alpha$ and $\beta_{ij}$ with the associated size of the equivalence class and the $P_4$ value, which is sometimes zero.
Table \ref{table:twosevensdelta} shows the duals of the corresponding non-equivalent values for $\epsilon$ and $\delta_{ij}$.
Finally, table \ref{table:twosevensequiv} shows the sizes of the combined $\beta-\delta$ equivalence classes.
These values allow to compute $D(4)$.
\end{example}

\begin{table}
   \center
\begin{tabular}{|c|c|c|c|c|c|c|c|c|c|}
\hline
     & $\alpha$ & $\beta_{12}$ & $\beta_{13}$ & $\beta_{14}$ & $\beta_{23}$ & $\beta_{24}$ & $\beta_{34}$ & \#EQ & P4 \\ \hline
A1 & $\bot$ & $\bot$ & $\bot$ & $\bot$ & $\bot$ & $\bot$ & $\bot$ & 1 & 1 \\ \hline
A2 & $\bot$ & ${0}$ & $\bot$ & $\bot$ & $\bot$ & $\bot$ & $\bot$ & 6 & 0 \\ \hline
A3 & $\bot$ & ${0}$ & ${0}$ & $\bot$ & $\bot$ & $\bot$ & $\bot$ & 15 & 0 \\ \hline
A4 & $\bot$ & ${0}$ & ${0}$ & ${0}$ & $\bot$ & $\bot$ & $\bot$ & 4 & 1 \\ \hline
A5 & $\bot$ & ${0}$ & ${0}$ & ${0}$ & ${0}$ & $\bot$ & $\bot$ & 15 & 0 \\ \hline
A6 & $\bot$ & ${0}$ & ${0}$ & ${0}$ & ${0}$ & ${0}$ & $\bot$ & 6 & 1 \\ \hline
A7 & $\bot$ & ${0}$ & ${0}$ & ${0}$ & ${0}$ & ${0}$ & ${0}$ & 1 & 4 \\ \hline
A8 & ${0}$ & ${0}$ & ${0}$ & ${0}$ & ${0}$ & ${0}$ & ${0}$ & 1 & 1 \\ \hline

\end{tabular}
\caption{Non equivalent values for $\beta_{ij}$ in the system of Theorem \ref{the:twotimesseven}}
\label{table:twosevensbeta}
\end{table}

\begin{table}
   \center
\begin{tabular}{|c|c|c|c|c|c|c|c|c|c|}
\hline

 & $\epsilon$ & $\delta_{12}$ & $\delta_{13}$ & $\delta_{14}$ & $\delta_{23}$ & $\delta_{24}$ & $\delta_{34}$ & \#EQ & P4 \\ \hline
B1 & ${0}$ & ${0}$ & ${0}$ & ${0}$ & ${0}$ & ${0}$ & ${0}$ & 1 & 1 \\ \hline
B2 & ${0}$ & ${0}$ & ${0}$ & ${0}$ & ${0}$ & ${0}$ & $\bot$ & 6 & 0 \\ \hline
B3 & ${0}$ & ${0}$ & ${0}$ & ${0}$ & ${0}$ & $\bot$ & $\bot$ & 15 & 0 \\ \hline
B4 & ${0}$ & ${0}$ & ${0}$ & ${0}$ & $\bot$ & $\bot$ & $\bot$ & 4 & 1 \\ \hline
B5 & ${0}$ & ${0}$ & ${0}$ & $\bot$ & $\bot$ & $\bot$ & $\bot$ & 15 & 0 \\ \hline
B6 & ${0}$ & ${0}$ & $\bot$ & $\bot$ & $\bot$ & $\bot$ & $\bot$ & 6 & 1 \\ \hline
B7 & ${0}$ & $\bot$ & $\bot$ & $\bot$ & $\bot$ & $\bot$ & $\bot$ & 1 & 4 \\ \hline
B8 & $\bot$ & $\bot$ & $\bot$ & $\bot$ & $\bot$ & $\bot$ & $\bot$ & 1 & 1 \\ \hline

\end{tabular}
\caption{Non equivalent values for $\widetilde\delta_{ij}$ in the system of Theorem \ref{the:twotimesseven}}
\label{table:twosevensdelta}
\end{table}

\begin{table}
   \center
\begin{tabular}{|c|c|c|c|c|c|c|c|c|}
\hline
\#EQ & B1 & B2 & B3 & B4 & B5 & B6 & B7 & B8 \\ \hline
A1   & 1 &  &  & 4 &  & 6 & 1 & 1 \\ \hline
A2   &  &  &  &  &  &  &  &  \\ \hline
A3   &  &  &  &  &  &  &  &  \\ \hline
A4   & 4 &  &  & 4 &  &  &  &  \\ \hline
A5   &  &  &  &  &  &  &  &  \\ \hline
A6   & 6 &  &  &  &  &  &  &  \\ \hline
A7   & 1 &  &  &  &  &  &  &  \\ \hline
A8   & 1 &  &  &  &  &  &  &  \\ \hline
\end{tabular}
\caption{Number of equivalent combinations, leading to a non-zero contribution, in the system of Theorem \ref{the:twotimesseven}}
\label{table:twosevensequiv}
\end{table}

Note the large number of empty cells, corresponding with terms not to be counted in Theorem \ref{the:twotimesseven}. Applying this theorem, the tables \ref{table:twosevensbeta}, \ref{table:twosevensdelta}, \ref{table:twosevensequiv}, leads to the following result, where the exponents in the powers of two correspond to the number of intervals $[\beta_{ij},\delta_{ij}]$ of size two:
\begin{gather*}
D(4) = 2^6 + 2 \times (4 \times 2^3 + 6 \times 2^1 + 4 \times 2^0 +  2^1 + 2 \times 2^0)\\
 = 2^6 + 8 \times 2^3 + 14 \times 2^1 + 12 \times 2^0 = 168
\end{gather*}

\section{Conclusion}
\label{sec:conclusions}
We presented a number of systems of equations in the free distributive lattice of antichains and introduced a technique to count the number of solutions.
The systems of equations were chosen to be applicable to the problem of counting the number of antichains of subsets of a finite set of $n$ elements, the $n^{th}$ Dedekind number $D(n)$ which is the size of the antichain set $\D(n)$.
We could apply our technique to reduce $D(n+2)$ to a sum over products of interval sizes in $\D(n)$, which was used in our computation of the ninth Dedekind number.
Finally, we showed that the technique can be generalised for similar reductions of $D(n + 3)$ and
$D(n + 4)$. 

In combination with the matrix computing formalism as developed by Christian J\"akel  \cite{JAKEL2023100006}, the techniques presented here may help to further reduce the complexity of the problem. Also FPGA implementations as we realised in \cite{10.1145/3674147}, could be developed to compute the P-coefficient independently. These and other innovative improvements, together with new developments in computing technology, may make it possible to compute $D(10)$ between now and 2050 as it has been happening all over the previous century.

Another question is the generalisation of the computational techniques very specifically developed to count monotone Boolean functions to other or more general lattices, see e.g. the recent publication of Berman and K\"ohler \cite{BERMANKOHLER2021}.


\bibliography{PcoefBIB}

\end{document}